\documentclass[cameraready]{jcmsi}%

\usepackage[dvips]{graphicx}
\usepackage{latexsym}
\usepackage{amssymb,amsmath}
\usepackage{url}

\usepackage{txfonts}
\makeatletter%
\input{ot1txtt.fd}
\makeatother%

\newcommand{\R}{\mathbb{R}}	
\newcommand{\C}{\mathbb{C}}

\newcommand{\N}{\mathbb{N}}
\newcommand{\be}{\begin{equation}}
\newcommand{\ee}{\end{equation}}	
\newcommand{\Pin}{\mathrm{Pin}}	
\newcommand{\Spin}{\mathrm{Spin}}
\renewcommand{\O}{\mathrm{O}}
\newcommand{\SO}{\mathrm{SO}}
\renewcommand{\eqref}[1]{(\ref{#1})}
\newcommand{\gvec}[1]{\ensuremath{\mbox{\textbf{\textit{#1}}}}}
\newcommand{\vect}[1]{\ensuremath{\mbox{\textbf{\textit{#1}}}}}

\Vol{4}
\No{1}
\Month{1}
\Year{2011}

\title{Introduction to Clifford's Geometric Algebra}

\begin{document}

\AUTHOR{%
  \author{Eckhard Hitzer}{aut:EH}
}

\AFFILIATE{%
  \affiliate{Department of Applied Physics, 
    University of Fukui, 
    Fukui 910-8507, Japan}{aut:EH}
}

\email{hitzer@mech.u-fukui.ac.jp}

\begin{abstract}%
Geometric algebra was initiated by W.K. Clifford over 130 years ago. It unifies all branches of physics, and has found rich applications in robotics, signal processing, ray tracing, virtual reality, computer vision, vector field processing, tracking, geographic information systems and neural computing. This tutorial explains the basics of geometric algebra, with concrete examples of the plane, of 3D space, of spacetime, and the popular conformal model. Geometric algebras are ideal to represent geometric transformations in the general framework of Clifford groups (also called versor or Lipschitz groups). Geometric (algebra based) calculus allows, e.g., to optimize learning algorithms of Clifford neurons, etc. 
\end{abstract}


\received{0}{00}{2011}
\revised{0}{00}{2011}

\maketitle


\section{Introduction}

W.K. Clifford (1845-1879), a young English Goldsmid professor of applied mathematics at the University College of London, published in 1878 in the American Journal of Mathematics Pure and Applied a nine page long paper on \textit{Applications of Grassmann's Extensive Algebra}. In this paper, the young genius Clifford, standing on the shoulders of two giants of algebra: W.R. Hamilton (1805-1865), the inventor of \textit{quaternions}, and H.G. Grassmann (1809-1877), the inventor of \textit{extensive algebra}, added the measurement of length and angle to Grassmann's abstract and coordinate free algebraic methods for computing with a space and all its subspaces. Clifford thus unified and generalized in his geometric algebras (=Clifford algebras) the works of Hamilton and Grassmann by finalizing the fundamental concept of \textit{directed numbers} \cite{GS:CMinGA}.

Any Clifford algebra $Cl(V)$ is generated from an inner-product\footnote{The inner product defines the measurement of length and angle.} vector space ($V$, $a\cdot b: a,b \in V \mapsto \R$)
by Clifford's geometric product setting\footnote{This setting amounts to an algebra generating relationship.} the geometric product\footnote{No product sign will be introduced, simple juxtaposition implies the geometric product just like $2x = 2 \times x$.} of any vector with itself equal to their inner product: $aa = a\cdot a$. We indeed have the \textit{universal property} \cite{HL:IAandGR,IP:CACG} that any isometry\footnote{
A $\mathbb{K}$-isometry between two inner-product spaces is a $\mathbb{K}$-linear mapping preserving the inner products. 
} 
from the vector space $V$ into an inner-product algebra\footnote{
A $\mathbb{K}$-algebra is a $\mathbb{K}$-vector space equipped with an associative and multilinear product.
An inner-product $\mathbb{K}$-algebra is a $\mathbb{K}$-algebra equipped with an inner product structure when taken as $\mathbb{K}$-vector space.
} 
$\mathcal{A}$ over the field\footnote{Important fields are real  $\mathbb{R}$ and complex numbers $\mathbb{C}$, etc. } $\mathbb{K}$ can be uniquely extended to an isometry\footnote{
That is a $\mathbb{K}$-linear homomorphism preserving the inner products, i.e., a $\mathbb{K}$-linear mapping
preserving both the products of the algebras when taken as rings,
and the inner products of the algebras when taken as inner-product vector spaces.} 
from the Clifford algebra $Cl(V)$ into $\mathcal{A}$. The Clifford algebra $Cl(V)$ is the \textit{unique} associative and multilinear algebra with this property. 
Thus if we wish to generalize methods from algebra, analysis, calculus, differential geometry (etc.) of real numbers, complex numbers, and quaternion algebra to vector spaces and multivector spaces (which include additional elements representing $2$D up to $n$D subspaces, i.e. plane elements up to hypervolume elements), the study of Clifford algebras becomes \textit{unavoidable}. Indeed repeatedly and independently a long list of Clifford algebras, their subalgebras and in Clifford algebras embedded algebras (like \textit{octonions} \cite{PL:CAandSpin}) of many spaces have been studied and applied historically, often under different names.

Some of these algebras are \textit{complex numbers (and the complex number plane), hyperbolic numbers (split complex numbers, real tessarines), dual numbers, quaternions, biquaternions (complex quaternions), dual quaternions, Pl\"{u}cker coordinates, bicomplex numbers (commutative quaternions, tessarines, Segre quaternions), Pauli algebra (space algebra), Dirac algebra (space-time algebra, Minkowski algebra), algebra of physical space, para-vector algebra, spinor algebra, Lie algebras, Cartan algebra, versor algebra, rotor algebra, motor algebra, Clifford bracket algebra, conformal algebra, algebra of differential forms, etc.}

Section \ref{ssc:defns} can also be skipped by readers less interested in mathematical definitions. Then follow sections on the geometric algebras of the plane $Cl(2,0)$, of 3D space $Cl(3,0)$, both with many notions also of importance in higher dimensions, and of spacetime $Cl(1,3)$, conformal geometric algebra $Cl(4,1)$, and finally on Clifford analysis. For further study we recommend e.g.~\cite{DFM:GAfCS}. 

\subsection{Definitions \label{ssc:defns}}

\textit{Definition of an algebra} \cite{Wiki:AoF}:
Let $\mathcal{A}$ be a vector space over the reals $\R$ with an additional binary operation from $\mathcal{A} \times \mathcal{A}$ to $\mathcal{A}$, denoted here by $\circ$ ($x \circ y$ is the product of any $x,y \in \mathcal{A}$). Then $\mathcal{A}$ is an algebra over $\R$ if the following identities hold $\forall \, x,y,z \in \mathcal{A}$, and ``scalars" $\alpha, \beta \in \R$ :
(1,2) Left and right distributivity: 
    $(x + y) \circ z = x \circ z + y \circ z$,
    $x \circ (y + z) = x \circ y + x \circ z$.
(3) Compatibility with scalars: 
    $(\alpha x) \circ (\beta y) = (\alpha\beta)(x \circ y)$.
This means that $x \circ y$ is bilinear. The binary operation is often referred to as multiplication in $\mathcal{A}$, which is not necessarily associative.

\textit{Definition of inner product space} \cite{Wiki:IP}:
An inner product space is a vector space $V$ over $\R$ together with an inner product map
$\langle . , . \rangle : V \times V \rightarrow \R$,
that satisfies $\forall \, x,y,z \in V$ and $\forall \, \alpha \in \R$:
(1) Symmetry:
    $\langle x,y \rangle = \langle y,x \rangle $.
(2) Linearity in the first argument:
    $\langle \alpha x,y\rangle= \alpha \langle x,y\rangle, \,\,\,
        \langle x+y,z\rangle= \langle x,z\rangle+ \langle y,z\rangle.$

\textit{Note:} We do not assume positive definiteness. 

\textit{Definition of inner product algebra}:
An inner product algebra, is an algebra equipped with an inner product $\mathcal{A} \times \mathcal{A} \rightarrow \R$.

\textit{Definition of Clifford's geometric algebra (GA)} \cite{FM:ICNAAM2007,PL:CAandSpin}: 
Let $\{e_1, e_2, \ldots , e_p, e_{p+1}, \ldots , e_{p+q}, e_{p+q+1}, \ldots, e_n \}$, with $n=p+q+r$, $e_k^2=\varepsilon_k$, $\varepsilon_k = +1$ for $k=1, \ldots , p$, $\varepsilon_k = -1$ for $k=p+1, \ldots , p+q$, $\varepsilon_k = 0$ for $k=p+q+1, \ldots , n$,  be an \textit{orthonormal base} of the inner product vector space $\R^{p,q,r}$ with a geometric product according to the multiplication rules 
\be
  e_k e_l + e_l e_k = 2 \varepsilon_k \delta_{k,l}, 
  \qquad k,l = 1, \ldots n,
\label{eq:mrules}
\ee
where $\delta_{k,l}$ is the Kronecker symbol with $\delta_{k,l}= 1$ for $k=l$, and $\delta_{k,l}= 0$ for $k\neq l$. This non-commutative product and the additional axiom of \textit{associativity} generate the $2^n$-dimensional Clifford geometric algebra $Cl(p,q,r) = Cl(\R^{p,q,r}) = Cl_{p,q,r} = \mathcal{G}_{p,q,r} = \R_{p,q,r}$ over $\R$. The set $\{ e_A: A\subseteq \{1, \ldots ,n\}\}$ with $e_A = e_{h_1}e_{h_2}\ldots e_{h_k}$, $1 \leq h_1< \ldots < h_k \leq n$, $e_{\emptyset}=1$, forms a graded (blade) basis of $Cl(p,q,r)$. The grades $k$ range from $0$ for scalars, $1$ for vectors, $2$ for bivectors, $s$ for $s$-vectors, up to $n$ for pseudoscalars. 
The vector space $\R^{p,q,r}$ is included in $Cl(p,q,r)$ as the subset of 1-vectors. The general elements of $Cl(p,q,r)$ are real linear combinations of basis blades $e_A$, called Clifford numbers, multivectors or hypercomplex numbers.

\textit{Note:} The definition of Clifford's GA is fundamentally \textit{coordinate system independent}, i.e. \textit{coordinate free}. Equation (\ref{eq:mrules}) is fully equivalent to the coordinate independent definition: $aa = a\cdot a, \, \forall \, a \in \R^{p,q,r}$ \cite{HS:CAtoGC}. This even applies to Clifford analysis (Geometric Calculus). Clifford algebra is thus ideal for computing with geometrical \textit{invariants} \cite{HL:IAandGR}. A review of five different ways to define GA, including one definition based on vector space basis element multiplication rules, and one definition focusing on $Cl(p,q,r)$ as a \textit{universal} associative algebra, is given in chapter 14 of the textbook \cite{PL:CAandSpin}.

In general $\langle A \rangle_{k}$ denotes the grade $k$ part of $A\in Cl(p,q,r)$. The parts of grade $0$, $(s-k)$, $(k-s)$, and $(s+k)$, respectively, of the geometric product of a $k$-vector $A_k\in Cl(p,q,r)$ with an $s$-vector $B_s\in Cl(p,q,r)$ 
\begin{gather}
  \label{eq:gaprods}
  A_k \ast B_s := \langle A_k B_s \rangle_{0}, \quad
  A_k \rfloor B_s := \langle A_k B_s \rangle_{s-k}, 
  \\
  A_k \lfloor B_s := \langle A_k B_s \rangle_{k-s}, \quad
  A_k \wedge B_s := \langle A_k B_s \rangle_{k+s},
\nonumber
\end{gather}
are called \textit{scalar product}, \textit{left contraction} (zero for $s<k$), \textit{right contraction} (zero for $k<s$), and (associative) \textit{outer product}, respectively. These definitions extend by linearity to the corresponding products of general multivectors. The various derived products of (\ref{eq:gaprods}) are related, e.g. by
\begin{align} 
  (A\wedge B)\rfloor C &= A\rfloor (B\rfloor C),
  \nonumber \\
  \forall\, A,B,C &\in Cl(p,q,r).
  \label{eq:ABCrel}
\end{align} 
Note that for vectors $a,b$ in $\R^{p,q,r} \subset Cl(p,q,r)$ we have
\begin{align} 
  ab &= a \rfloor b + a \wedge b, 
  \nonumber \\
  a \rfloor b &= a \lfloor b = a \cdot b = a \ast b,
\end{align} 
where $a \cdot b$ is the usual inner product of $\R^{p,q,r}$. 

For $r=0$ we often denote $\R^{p,q}=\R^{p,q,0}$, and $Cl(p,q) = Cl(p,q,0)$. For Euclidean vector spaces $(n=p)$ we use $\R^{n}=\R^{n,0}=\R^{n,0,0}$, and $Cl(n) = Cl(n,0) = Cl(n,0,0)$. The \textit{even} grade subalgebra of $Cl(p,q,r)$ is denoted by $Cl^+(p,q,r)$, the $k$-vectors of its basis have only even grades $k$. Every $k$-vector $B$ that can be written as the outer product $B = b_1 \wedge b_2 \wedge \ldots \wedge b_k$ of $k$ vectors $b_1, b_2, \ldots, b_k \in \R^{p,q,r}$ is called a \textit{simple} $k$-vector or \textit{blade}. 

\textit{Definition of outermorphism} \cite{AM:LinGA,HS:CAtoGC}: 
An \textit{outermorphism} is the unique extension to $Cl(V)$ of a vector space map for all $a\in V$, $f: a\mapsto f(a) \in V'$, and is given by the mapping $B = b_1 \wedge b_2 \wedge \ldots \wedge b_k \rightarrow f(b_1) \wedge f(b_2) \wedge \ldots \wedge f(b_k)$ in $Cl(V')$, for every blade $B$ in $Cl(V)$.

\section{Clifford's GA of a Euclidean plane}

In order to demonstrate how to compute with Clifford numbers, we begin with a low dimensional example. 

\subsection{Example of $Cl(2,0)$}
A Euclidean plane $\R^2=\R^{2,0}=\R^{2,0,0}$ is spanned by $e_1, e_2 \in \R^2$ with 
\begin{gather}
  e_1\cdot e_1 = e_2\cdot e_2 = 1, \quad e_1\cdot e_2 = 0.
\end{gather}
$\{e_1, e_2\}$ is an \textit{orthonormal} vector basis of $\R^2$. 

Under Clifford's \textit{associative} geometric product we set
\begin{gather}
  e_1^2=e_1e_1 := e_1\cdot e_1=1, 
  \nonumber \\ 
  e_2^2=e_2e_2 := e_2\cdot e_2=1, 
  \\
  \hspace*{-12mm}\text{and}\quad 
  (e_1+e_2)(e_1+e_2) 
  = e_1^2+e_2^2+e_1e_2 + e_2e_1 \\
  = 2 + e_1e_2 + e_2e_1
  := (e_1+e_2)\cdot(e_1+e_2) = 2.
  \nonumber 
\end{gather}
Therefore
\be 
  e_1e_2 + e_2e_1 = 0 
  \,\,\,\Leftrightarrow \,\,\,
  e_1e_2 = - e_2e_1,
\ee 
i.e. the geometric product of orthogonal vectors forms a new entity, called unit \textit{bi-vector} $e_{12}=e_1e_2$ by Grassmann, and is \textit{anti-symmetric}. General bivectors in $Cl(2,0)$ are $\beta e_{12}$, $\forall \, \beta \in \R \setminus \{0\}$. For orthogonal vectors the geometric product equals Grassmann's anti-symmetric outer product (exterior product, symbol $\wedge$) 
\begin{align} 
  e_{12} &= e_1e_2 = e_1 \wedge e_2 
  \nonumber \\
  &= - e_2 \wedge e_1 = -e_2e_1 = -e_{21}.
\end{align} 
Using associativity, we can compute the products
\be   
  e_1 e_{12} = e_1 e_1 e_2 = e_1^2 e_2 = e_2, \quad
  e_2 e_{12} = - e_2 e_{21} = - e_1,
\ee 
which represent a mathematically \textit{positive} (anti-clockwise) $90^{\circ}$ \textit{rotation}. The opposite order gives
\be 
  e_{12}e_1 = -e_{21}e_1 = -e_2, \quad e_{12}e_2 = e_1,
\ee 
which represents a mathematically \textit{negative} (clockwise) $90^{\circ}$ \textit{rotation}. The bivector $e_{12}$ acts like a \textit{rotation operator}, and we observe the general anti-commutation property
\be 
  a e_{12} = - e_{12} a, \quad \forall a=a_1e_1+a_2e_2 \in \R^2, \,\,\,a_1,a_2 \in \R.
\ee 
The square of the unit bivector is $-1$,
\be 
  e_{12}^2 = e_1 e_2 e_{12} = e_1 (-e_1) = -1,
\ee 
just like the imaginary unit $j$ of complex numbers $\C$. 

Table \ref{tb:Cl2mtable} is the complete multiplication table of the Clifford algebra $Cl(\R^2) = Cl(2,0,0) = Cl(2,0)$ with algebra basis elements $\{1, e_1, e_2, e_{12}\}$ (which includes the vector basis of $\R^2$). The even subalgebra spanned by $\{1,e_{12}\}$ (closed under geometric multiplication), consisting of even grade scalars (0-vectors) and bivectors (2-vectors), is isomorphic to $\mathbb{C}$.

\subsection{Algebraic unification and vector inverse}

The general geometric product of two vectors $a,b \in \R^2$ 
\begin{gather}
  ab=(a_1e_1+a_2e_2)(b_1e_1+b_2e_2) \nonumber \\
  = a_1b_1 + a_2b_2 + (a_1b_2 - a_2b_1) e_{12} \\
  = \frac{1}{2}(ab+ba) + \frac{1}{2}(ab-ba)
  = a \cdot b + a\wedge b,
  \nonumber 
\end{gather} 
has therefore a scalar \textit{symmetric} inner product part
\begin{align} 
  \frac{1}{2}(ab+ba) 
  &= a \cdot b  
  =a_1b_1 + a_2b_2 
  \nonumber \\
  &= |a| |b| \cos \theta_{a,b},
\end{align} 
and a bi-vector \textit{skew-symmetric} outer product part
\be 
  \frac{1}{2}(ab-ba)=a\wedge b 
  =(a_1b_2 - a_2b_1) e_{12} = |a| |b| e_{12}\sin \theta_{a,b}. 
\ee
We observe that parallel vectors $(\theta_{a,b}=0)$ commute, $ab = a\cdot b = ba$, and orthogonal vectors $(\theta_{a,b}=90^{\circ})$ anti-commute, $ab = a \wedge b = -ba$. 
The outer product part $a\wedge b$ represents the \textit{oriented area} of the parallelogram spanned by the vectors $a,b$ in the plane of $\R^2$, with oriented magnitude 
\be 
\det(a,b) = |a| |b| \sin \theta_{a,b} = (a\wedge b) e_{12}^{-1},
\ee
where $e_{12}^{-1} = -e_{12}$, because $e_{12}^2=-1$.  

\begin{table}
\caption{Multiplication table of plane Clifford algebra $Cl(2,0)$.}
\label{tb:Cl2mtable}
\begin{center}
\begin{tabular}{c|cccc}
       & 1 & $e_1$ & $e_2$ & $e_{12}$\\
 \hline
 1     & 1 & $e_1$ & $e_2$ & $e_{12}$\\
 $e_1$ & $e_1$ & 1 & $e_{12}$ & $e_{2}$\\
 $e_2$ & $e_2$ & $-e_{12}$ & 1 & $-e_{1}$\\
 $e_{12}$ & $e_{12}$ & $-e_{2}$ & $e_1$ & $-1$ 
\end{tabular}
\end{center}
\end{table}

With the \textit{Euler} formula we can rewrite the geometric product as
\begin{align} 
  ab &= |a| |b| (\cos \theta_{a,b} + e_{12}\sin \theta_{a,b}) 
  \nonumber \\
  &= |a| |b| e^{\theta_{a,b} e_{12}} \,,
  \label{eq:Euler}
\end{align} 
again because $e_{12}^2=-1$. 

The geometric product of vectors is \textit{invertible} for all vectors with non-zero square 
$a^2 \neq 0$
\begin{gather}
  a^{-1} := a/a^2, \quad 
  a a^{-1} = a a/a^2 = 1, \nonumber \\
  a^{-1} a = \frac{a}{a^2}a = a^2/a^2 = 1. 
\end{gather} 
The inverse vector $a/a^2$ is a rescaled version (reflected at the unit circle) of the vector $a$. 
This invertibility leads to enormous simplifications and ease of practical computations.

\subsection{Geometric operations and transformations}

For example, the \textit{projection} of one vector $x\in\R^2$ onto another $a \in \R^2$ is
\be 
  x_{\parallel} = | x \,| \cos\theta_{a,x}\frac{a}{|a|} 
  = (x \cdot \frac{a}{|a|})\frac{a}{|a|}
  = (x \cdot a)\frac{a}{|a|^2} 
  = (x \cdot a)a^{-1}.
\ee 
The \textit{rejection} (perpendicular part) is
\begin{gather} 
  x_{\perp} = x-x_{\parallel}
  = x aa^{-1} - (x \cdot a)a^{-1} 
  \nonumber \\
  = (xa-x \cdot a)a^{-1}
  = (x\wedge a)a^{-1}.
\end{gather}

We can now use $x_{\parallel}, x_{\perp}$ to compute the reflection\footnote{Note that reflections at hyperplanes are nothing but the \textit{Householder transformations} \cite{Wiki:Householder} of matrix analysis.} of $x=x_{\parallel}+ x_{\perp}$ at the line (hyperplane\footnote{A hyperplane of a $n$D space is a $(n-1)$D subspace, thus a hyperplane of $\R^2$, $n=2$, is a 1D ($2-1=1$) subspace, i.e. a line. Every hyperplane is characterized by a vector normal to the hyperplane.}) with normal vector $a$, which means to reverse $x_{\parallel}\rightarrow -x_{\parallel}$
\begin{gather}
  x' = -x_{\parallel} + x_{\perp} 
  = -a^{-1}a\,x_{\parallel} + a^{-1}a\,x_{\perp} 
  \nonumber \\
  = -a^{-1}x_{\parallel}a - a^{-1}\,x_{\perp}a
  = -a^{-1}(x_{\parallel} +\,x_{\perp})a
  = -a^{-1} x a .
  \label{eq:refl}
\end{gather}
The combination of two reflections at two lines (hyperplanes) with normals $a,b$
\be 
  x'' = -b^{-1} x' b = b^{-1} a^{-1} x ab = (ab)^{-1} x ab = R^{-1} x R ,
\ee 
gives a rotation. The rotation angle is $\alpha = 2 \theta_{a,b}$ and the \textit{rotor}
\be
  R = e^{\theta_{a,b} e_{12}} = e^{\frac{1}{2}\alpha e_{12}},
\ee 
where the lengths $|a||b|$ of $ab$ cancel against 
$|a|^{-1}|b|^{-1}$ in $(ab)^{-1}$.
The rotor $R$ gives the \textit{spinor} form of rotations, fully replacing rotation matrices, and introducing the same elegance to \textit{real} rotations in $\R^2$, like in the complex plane. 

In 2D, the product of three reflections, i.e. of a rotation and a reflection, leads to another reflection. In 2D the product of an \textit{odd} number of reflections always results in a \textit{reflection}. That the product of an \textit{even} number of reflections leads to a \textit{rotation} is true in general dimensions. These transformations are in Clifford algebra simply described by the products of the vectors normal to the lines (hyperplanes) of reflection and called versors. 

\textit{Definition of a versor} \cite{HL:IAandGR}:
A \textit{versor} refers to a Clifford monomial (product expression) composed of invertible vectors. It is called a \textit{rotor}, or \textit{spinor}, if the number of vectors is even. It is called a \textit{unit versor} if its magnitude is 1. 

Every versor $A = a_1\ldots a_r, \quad a_1, \ldots, a_r \in \R^2, r\in \N$ has an inverse 
\be 
  A^{-1} = a_r^{-1}\ldots a_1^{-1} = a_r\ldots a_1/(a_1^2 \ldots a_r^2), 
\ee 
such that
\be 
  A A^{-1} = A^{-1} A = 1.
\ee 
This makes the set of all versors in $Cl(2,0)$ a group, the so called \textit{Lipschitz group} with symbol $\Gamma(2,0)$, also called \textit{Clifford group} or \textit{versor group}. Versor transformations apply via \textit{outermorphisms} to all elements of a Clifford algebra. It is the group of all reflections and rotations of $\R^2$. The reverse product order of a versor represents an involution (applying it twice leads to identity) called \textit{reversion}\footnote{Reversion is an anti-automorphism. Often a dagger $A^{\dagger}$ is used instead of the tilde, as well as the term transpose.}
\be 
  \widetilde{A} = (a_1\ldots a_r)^{\sim} = a_r\ldots a_1. 
\ee 
In the case of $Cl(2,0)$ we have 
\begin{gather}
  \Gamma(2,0) = Cl^-(2,0) \cup Cl^+(2,0)
\end{gather}
where the odd grade vector part $Cl^-(2,0)=\R^2$ generates reflections, and the even grade part of scalars and bivectors 
$Cl^+(2,0)= \{A \mid A = \alpha + \beta e_{12}, \alpha,\beta \in \R \}$
 generates rotations. 
The normalized subgroup of versors is called \textit{pin group}
\be 
  \Pin(2,0) = \{A\in \Gamma(2,0) \mid A\widetilde{A}= \pm 1\}.
\ee 
In the case of $Cl(2,0)$ we have 
\begin{align}
  &\Pin(2,0) 
  \\
  &= \{a \in \R^2 \mid a^2=1 \} \cup 
  \{A \mid A = \cos \varphi + e_{12} \sin \varphi , \varphi \in \R \}.
    \nonumber 
\end{align} 
The pin group has an even subgroup, called \textit{spin group}
\be 
  \Spin(2,0) = \Pin(2,0) \cap Cl^+(2,0).
\ee 
As mentioned above $Cl^+(2,0)$ has the basis $\{1, e_{12}\}$ and is thus isomorphic to $\C$. 
In the case of $Cl(2,0)$ we have explicitly
\begin{align} 
  &\Spin(2,0) 
  \\
  &= \{A \mid  A = \cos \varphi + e_{12} \sin \varphi , \varphi \in \R \}.
  \nonumber 
\end{align}
The spin group has in general a \textit{spin plus subgroup}\footnote{Note, that in general for Clifford algebras $Cl(n,0)$ of Euclidean spaces $\R^{n,0}$ we have the identity 
$\Spin(n) = \Spin_+(n)$, where $\Spin(n) = \Spin(n,0)$. The reason is that $A\widetilde{A} < 0$ is only possible for non-Euclidean spaces $\R^{p,q}$, with $q>0$.}
\be 
  \Spin_+(2,0) = \{A\in \Spin(2,0) \mid A\widetilde{A}= +1\}. 
\ee
The groups $\Pin(2,0)$, $\Spin(2,0)$, and $\Spin_+(2,0)$ are two-fold coverings\footnote{Two-fold covering means, that there are always two elements $\pm A$ in $\Pin(2,0)$, $\Spin(2,0)$, and $\Spin_+(2,0)$, representing one element in $\O(2,0)$, $\SO(2,0)$ and $\SO_+(2,0)$, respectively.} of the orthogonal group $\O(2,0)$ \cite{Wiki:IOG}, the special orthogonal group $\SO(2,0)$, and the component of the special orthogonal group connected to the identity $\SO_+(2,0)$. 

Let us point out, that this natural combination of reflections leading to spinors rightly indicates the way how Clifford's GA is able to give a fully real algebraic description to quantum mechanics, with a clear cut geometric interpretation.

\subsection{Vectors, $k$-vectors and multivectors}

A general element in $Cl(2,0)$, also called \textit{multivector} can be represented as
\begin{align} 
  M =\, &m_0 + m_1 e_1 + m_2 e_2 + m_{12} e_{12}, 
  \nonumber \\
  &m_0, m_1, m_2, m_{12} \in \R. 
\end{align}
Like real and imaginary parts of a complex number, we have 
a scalar part $\langle M \rangle_0$ of grade $0$, 
a vector part $\langle M \rangle_1$ of grade $1$, 
and a bivector part $\langle M \rangle_2$ of grade $2$
\be 
  M = \langle M \rangle_0 + \langle M \rangle_1 + \langle M \rangle_2,
\ee 
with 
\begin{align} 
  \langle M \rangle_0 &= m_0, \qquad
  \langle M \rangle_1 = m_1 e_1 + m_2 e_2, 
  \nonumber \\
  \langle M \rangle_2 &= m_{12} e_{12} \,. 
\end{align}
The set of all grade $k$ elements, $0 \leq k \leq 2$, is denoted $Cl^k(2,0)$.

\textit{Grade extraction} $\langle \ldots \rangle_k$, $0 \leq k \leq 2$, allows to do many useful computations, like e.g. the \textit{angle} between two vectors (atan = $\tan^{-1}$ = arctan)
\be 
  \theta_{a,b} = \mathrm{atan} \frac{\langle ab \rangle_2 e_{12}^{-1}}{\langle ab \rangle_0}
               \stackrel{(\ref{eq:Euler})}{=} \mathrm{atan} \frac{\sin \theta_{a,b}}{\cos \theta_{a,b}} . 
\ee 
The symmetric scalar part of the geometric product of two multivectors $M,N \in Cl(2,0)$ is also called their \textit{scalar product} and because of its fundamental importance denoted with a special product sign $\ast$
\begin{align} 
  &M \ast N = \langle MN \rangle_0 
  \nonumber \\
  &= m_0n_0 + m_1n_1 + m_2 n_2 - m_{12}n_{12}=N\ast M ,
\end{align} 
which is easy to compute using the multiplication table Table \ref{tb:Cl2mtable}, since only the \textit{diagonal} entries of Table \ref{tb:Cl2mtable} are scalar. 
Using the reversion ($\widetilde{e_{12}}=\widetilde{e_{1}e_2}=e_2e_1 = -e_1e_2=-e_{12}$)
\be 
  \widetilde{M} = m_0 + m_1 e_1 + m_2 e_2 - m_{12} e_{12},
\ee  
we get the \textit{norm} $|M|$ of a multivector as 
\be 
  |M|^2 = M \ast \widetilde{M} = \langle M\widetilde{M} \rangle_0
  = m_0^2 + m_1^2 + m_2^2 + m_{12}^2 ,
  \label{eq:norm}
\ee 
which for vectors $M \in Cl^1(2,0)$ is identical to the length of a vector, and for even grade subalgebra elements $M \in Cl^+(2,0)$ to the modulus of complex numbers, and for pure bivectors like $a\wedge b = \langle ab \rangle_2$ gives the area content of the parallelogram spanned by the vectors $a,b \in \R^2$
\be 
  |a\wedge b| = |a|\,|b|\,|\sin \theta_{a,b}| = |\det(a,b)|.
\ee 

The reversion operation maps $\widetilde{e_{12}}=-e_{12}$, it is therefore the equivalent of complex conjugation in the isomorphism $Cl^+(2,0) \cong \C$, fully consistent with its use in the norm (\ref{eq:norm}).

\subsection{Generalizing inner and outer products to multivectors}

The grade extraction also allows us to generalize the inner product of vectors, which maps two grade one vectors $a,b \in \R^2=Cl^1(2,0)$ to a grade zero scalar
\be 
  a \cdot b = \langle ab \rangle_{(1-1=0)} \in \R=Cl^0(2,0),
\ee 
and therefore \textit{lowers} the grade by $1$. In contrast, the outer product with a vector \textit{raises} the grade by $1$, i.e.
\be 
  a \wedge b = \langle ab \rangle_{(1+1=2)} \in Cl^2(2,0).
\ee 

In general the \textit{left contraction} (symbol $\rfloor$) of a $k$-vector $A_k = \langle A \rangle_k$ with a $l$-vector 
$B_l = \langle B \rangle_l$ is defined as 
\be 
  A_k \rfloor B_l = \langle A_k B_l \rangle_{(l-k)} , 
\ee 
which is zero if $0>l-k$, i.e. if $k > l$. Figuratively speaking (like with projections), we can only contract objects of the same or lower dimension from the left onto an object on the right. 
The \textit{right contraction} (symbol $\lfloor$) is defined as
\be 
  A_k \lfloor B_l = \langle A_k B_l \rangle_{(k-l)} , 
\ee 
which is zero if $k-l < 0$, i.e. if $k < l$. Figuratively speaking, we can only contract objects of the same or lower dimension from the right onto an object on the left. Both contractions are bilinear and can thus be extended to contractions of multivectors
\begin{align}
  A \rfloor B &= \sum_{k=0}^2\sum_{l=0}^2 \langle \langle A \rangle_k \langle B \rangle_l \rangle_{(l-k)}, 
  \\
  A \lfloor B &= \sum_{k=0}^2\sum_{l=0}^2 \langle \langle A \rangle_k \langle B \rangle_l \rangle_{(k-l)}.
\end{align}
The reversion changes the order of factors and therefore relates left and right contractions by
\be 
  (A \rfloor B)^{\sim} = \widetilde{B} \,\lfloor \,\widetilde{A}, \qquad 
  (A \lfloor B)^{\sim} = \widetilde{B} \,\rfloor \,\widetilde{A}. 
\ee 

In general the associative outer product of a $k$-vector $A_k = \langle A \rangle_k$ with a $l$-vector $B_l = \langle B \rangle_l$ is defined as the maximum grade part of the geometric product
\be 
  A_k \wedge B_l = \langle A_k B_l \rangle_{(l+k)} , 
\ee 
where $A_k \wedge B_l=0$ for $l+k > 2$, because in $Cl(2,0)$ the highest grade possible is $2$. For example the following outer products are obtained $(\alpha, \beta \in \R, a,b,c \in \R^2)$
\begin{gather} 
  \alpha \wedge \beta = \alpha \beta, \qquad
  \alpha \wedge a = a \wedge \alpha = \alpha a, 
  \nonumber \\
  a \wedge b \wedge c = 0,  
\end{gather}  
where the last identity is specific to $Cl(2,0)$, it does generally not hold in GAs of higher dimensional vector spaces. The last identity is due to the fact, that in a plane every third vector $c$ can be expressed by linear combination of two linearly independent vectors $a,b$. If $a,b$ would not be linearly independent, then already $a\wedge b=0$, i.e. $a \parallel b$. The outer product in 2D is of great advantage, because the cross product of vectors does only exist in 3D, not in 2D. We next treat the GA of $\R^3$.

\section{Geometric algebra of 3D Euclidean space}

The Clifford algebra $Cl(\R^3) = Cl(3,0)$ is arguably the by far most thoroughly studied and applied GA. In physics it is also known as \textit{Pauli algebra}, since Pauli's spin matrices provide a $2\times 2$ matrix representation. This shows how GA unifies \textit{classical} with \textit{quantum} mechanics. 

Given an orthonormal vector basis $\{e_1, e_2, e_3\}$ of $\R^3$, the eight-dimensional ($2^3=8$) Clifford algebra $Cl(\R^3) = Cl(3,0)$ has a basis of one scalar, three vectors, three bivectors and one trivector
\be 
  \{1, e_1, e_2, e_3, e_{23}, e_{31}, e_{12}, e_{123}\},
  \label{eq:Cl3basis}
\ee 
where as before $e_{23}=e_2e_3, e_{123} = e_1e_2e_3$, etc. All basis bivectors square to $-1$, and the product of two basis bivectors gives the third
\be 
  e_{23} e_{31} = e_{21} = - e_{12}, \,\,\, \text{etc.}
\ee 
Therefore the even subalgebra $Cl^+(3,0)$ with basis\footnote{The minus signs are only chosen, to make the product of two bivectors identical to the third, and not minus the third.} $\{1, -e_{23}, -e_{31}, -e_{12}\}$ is indeed found to be isomorphic to quaternions $\{1, \mathbf{i}, \mathbf{j}, \mathbf{k}\}$. This isomorphism is not incidental. As we have learned already for $Cl(2,0)$, also in $Cl(3,0)$, the even subalgebra is the algebra of rotors (rotation operators) or spinors, and describes rotations in the same efficient way as do quaternions \cite{Wiki:EucGroup}. We therefore gain a \textit{real geometric} interpretation of quaternions, as the oriented bi-vector side faces of a unit cube, with edge vectors $\{e_1, e_2, e_3\}$. 

In $Cl(3,0)$ a reflection at a plane (=hyperplane) is specified by the plane's normal vector $a\in\R^3$
\be 
  x' = -a^{-1} x a,
\ee 
the proof is identical to the one in (\ref{eq:refl}) for $Cl(2,0)$. The combination of two such reflections leads to a rotation by $\alpha = 2 \theta_{a,b}$
\begin{align} 
  x'' &= R^{-1} x R, 
  \\ 
  R &= ab = |a||b| e^{\theta_{a,b}\mathbf{i}_{a,b}}
  = |a||b| e^{\frac{1}{2}\alpha \mathbf{i}_{a,b}},
  \nonumber
\end{align} 
where $\mathbf{i}_{a,b}={a \wedge b}/({|a\wedge b|})$ specifies the oriented unit bivector of the plane spanned by $a,b \in \R^3$. 

The unit trivector $i_3=e_{123}$ also squares to $-1$
\begin{align}
  i_3^2 &= e_1e_2e_3e_1e_2e_3 
  = - e_1e_2e_1e_3e_2e_3 \nonumber \\
  &= e_1e_2e_1e_2 e_3e_3 = (e_1e_2)^2 (e_3)^2 = -1,
\end{align}
where we only used that the permutation of two orthogonal vectors in the geometric product produces a minus sign. Hence $i_3^{-1} = - i_3$. We further find, that $i_3$ commutes with every vector, e.g.
\begin{align} 
  e_1 i_3 &= e_1 e_1e_2e_3 = e_{23}, 
  \\
  i_3 e_1 &= e_1e_2e_3 e_1 = - e_1e_2e_1e_3 = e_1e_1 e_2e_3 = e_{23}, 
  \nonumber 
\end{align} 
and the like for $e_2 i_3 = i_3 e_2$, $e_3 i_3 = i_3 e_3$. If $i_3$ commutes with every vector, it also commutes with every bivector $a\wedge b = \frac{1}{2}(ab-ba)$, hence $i_3$ commutes with every element of $Cl(3,0)$, a property which is called \textit{central} in mathematics. The central subalgebra spanned by $\{1,i_3\}$ is isomorphic to complex numbers $\C$. $i_3$ changes bivectors into orthogonal vectors
\be 
  e_{23} i_3 = e_2e_3e_1e_2e_3 = e_1 e_{23}^2 = -e_1\,, \,\,\, \text{etc.}
\ee 
This is why the GA $Cl(3,0)$ is isomorphic to \textit{complex quaternions}, one form of \textit{biquaternions}
\begin{gather} 
  \left\{ 1, e_1 = e_{23}i_3^{-1}, e_2 = e_{31}i_3^{-1}, e_3 = e_{12}i_3^{-1}, \right.
  \nonumber \\
  \left. e_{23}, e_{31}, e_{12}, i_3 \right\} \,. 
\end{gather} 
Yet writing the basis in the simple product form (\ref{eq:Cl3basis}), fully preserves the \textit{geometric interpretation} in terms of scalars, vectors, bivectors and trivectors, and allows to \textit{reduce} all products to elementary geometric products of basis vectors, which is used in \textit{computational optimization} schemes for GA software like Galoop \cite{HPK:Gaalop2010}. 

\begin{table}
\caption{Multiplication table of Clifford algebra $Cl(3,0)$ of Euclidean 3D space $\R^3$.}
\label{tb:Cl3mtable}
\begin{center}
\begin{tabular}{c|cccccccc}
       & $1$ & $e_1$ & $e_2$ & $e_3$ & $e_{23}$ & $e_{31}$ & $e_{12}$ & $e_{123}$\\
 \hline
 $1$     & $1$ & $e_1$ & $e_2$ & $e_3$ & $e_{23}$ & $e_{31}$ & $e_{12}$ & $e_{123}$\\
 $e_1$ & $e_1$ & 1 & $e_{12}$ & $-e_{31}$ & $e_{123}$ & $-e_3$ & $e_{2}$ & $e_{23}$\\
 $e_2$ & $e_2$ & $-e_{12}$ & 1 & $e_{23}$ & $e_3$ & $e_{123}$ & $-e_1$ & $e_{31}$ \\
 $e_3$ & $e_3$ & $e_{31}$ & $-e_{23}$ & $1$ & $-e_2$ & $e_1$ & $e_{123}$ & $e_{12}$\\
 $e_{23}$ & $e_{23}$ & $e_{123}$ & $-e_3$ & $e_2$ & $-1$ & $-e_{12}$ & $e_{31}$ & $-e_1$\\
 $e_{31}$ & $e_{31}$ & $e_3$ & $e_{123}$ & $-e_1$ & $e_{12}$ & $-1$ & $-e_{23}$ & $-e_2$\\
 $e_{12}$ & $e_{12}$ & $-e_{2}$ & $e_1$ & $e_{123}$ & $-e_{31}$ & $e_{23}$ & $-1$ & $-e_3$\\
 $e_{123}$ & $e_{123}$ & $e_{23}$ & $e_{31}$ & $e_{12}$ & $-e_1$ & $-e_2$ & $-e_3$ & $-1$
\end{tabular}
\end{center}
\end{table}

\subsection{Multiplication table and subalgebras of $Cl(3,0)$ }

For the full multiplication table of $Cl(3,0)$ we still need the geometric products of vectors and bivectors. By changing labels in Table \ref{tb:Cl2mtable} ($1 \leftrightarrow 3$ or $2 \leftrightarrow 3$), we get that
\begin{align} 
  e_2 e_{23} &= - e_{23}e_2 = e_3, 
  \nonumber \\
  e_3 e_{23} &= - e_{23}e_3 = -e_2
  \\
  e_1 e_{31} &= - e_{31}e_1 = -e_3, 
  \nonumber \\
  e_3 e_{31} &= - e_{31}e_3 = e_1,
\end{align} 
which shows that in general a vector and a bivector, which includes the vector, anti-commute. 
The products of a vector with its orthogonal bivector always gives the trivector $i_3$
\begin{gather} 
  e_1 e_{23} = e_{23}e_1 = i_3, \quad
  e_2 e_{31} = e_{31}e_2 = i_3,
  \nonumber \\
  e_3 e_{12} = e_{12}e_3 = i_3,
\end{gather} 
which also shows that in general vectors and orthogonal bivectors necessarily commute. 
Commutation relationships therefore clearly depend on both \textit{orthogonality} properties and on the \textit{grades} of the factors, which can frequently be exploited for computations even without the explicit use of coordinates. 

Table \ref{tb:Cl3mtable} gives the \textit{multiplication table} of $Cl(3,0)$. The elements on the left most column are to be multiplied from the left in the geometric product with the elements in the top row. 
Every subtable of Table \ref{tb:Cl3mtable}, that is closed under the geometric product represents a \textit{subalgebra} of $Cl(3,0)$. 

We naturally find Table \ref{tb:Cl2mtable} as a subtable, because with 
$\R^2 \subset \R^3$ we necessarily have $Cl(2,0) \subset Cl(3,0)$. In general any pair of orthonormal vectors will generate a $4$D \textit{subalgebra} of $Cl(3,0)$ \textit{specific to the plane} spanned by the two vectors.

We also recognize the subtable of the \textit{even subalgebra} with basis $\{1, e_{23}, e_{31}, e_{12}\}$ isomorphic to \textit{quaternions}. 

Then there is the subtable of $\{1, e_{123}\}$ of the \textit{central subalgebra} isomorphic to $\C$. Any element of $Cl(3,0)$ squaring to $-1$ generates a 2D subalgebra isomorphic to $\C$. This fundamental observation leads to the generalization of the conventional complex Fourier transformation (FT) to Clifford FTs, sometimes called geometric algebra FTs \cite{EH:QFTgen}. Moreover, we also obtain generalizations of real and complex (dual) wavelets to Clifford wavelets, which include quaternion wavelets.

We further have closed subtables with elements $\{1, e_1, e_{23}, e_{123}=i_3\}$, which are fully commutative and correspond to \textit{tessarines}, also called \textit{bicomplex numbers}, \textit{Segre quaternions}, \textit{commutative quaternions}, or \textit{Cartan subalgebras}. In general we can specify any unit vector $u\in\R^3$, $u^2=1$, and get a commutative tessarine subalgebra with basis $\{1, u, ui_3, i_3\}$, which includes besides $1$ and the vector $u$ itself, the bivector $ui_3$ orthogonal to $u$, and the oriented unit volume trivector $i_3$. Again GA provides a clear and useful geometric interpretation of tessarines and their products. Knowledge of the geometric interpretation is fundamental in order to \textit{see} these algebras in nature and in problems at hand to be solved, and to \textit{identify} settings, where a restriction to subalgebra computations may save considerable computation costs.

\subsection{The grade structure of $Cl(3,0)$ and duality}

A general multivector in $Cl(3,0)$, can be represented as
\begin{gather}
  M = m_0 + m_1 e_1 + m_2 e_2 + m_3 e_3 + m_{23} e_{23}
  + m_{31} e_{31}+ m_{12} e_{12}   \nonumber \\
  + m_{123} e_{123}, \quad
  m_0, \ldots , m_{123} \in \R. 
\end{gather}
We have 
a scalar part $\langle M \rangle_0$ of grade $0$, 
a vector part $\langle M \rangle_1$ of grade $1$, 
a bivector part $\langle M \rangle_2$ of grade $2$,
and a trivector part $\langle M \rangle_3$ of grade $3$
\begin{align}
  M &= \langle M \rangle_0 + \langle M \rangle_1 + \langle M \rangle_2 + \langle M \rangle_3,
  \\
   \langle M \rangle_0 &= m_0, \quad
  \langle M \rangle_1 = m_1 e_1 + m_2 e_2 + m_3 e_3, \nonumber \\
  \langle M \rangle_2 &= m_{23} e_{23} + m_{31} e_{31}+ m_{12} e_{12},   
  \quad
  \langle M \rangle_3 = m_{123} e_{123}.
  \nonumber 
\end{align}
The set of all grade $k$ elements, $0 \leq k \leq 3$, is denoted $Cl^k(3,0)$.

The multiplication table of $Cl(3,0)$, Table \ref{tb:Cl3mtable}, reveals that multiplication with $i_3$ (or $i_3^{-1}=-i_3$) consistently changes an element of grade $k$, $0 \leq k \leq 3$, into an element of grade $3-k$, i.e. scalars to trivectors (also called pseudoscalars) and vectors to bivectors, and vice versa. This means that the geometric product of a multivector $M \in Cl(3,0)$ and the pseudoscalar $i_3$ (or $i_3^{-1}=-i_3$) always results\footnote{In the context of blade subspaces, whenever a blade $B$ contains another blade $A$ as factor, then the geometric product is reduced to left or right contraction: $A B = A\rfloor B, 
  B A = B \lfloor A$.} in the left or right contraction
\be 
  M i_3 = M\rfloor i_3, 
  \qquad
  i_3 M = i_3 \lfloor M.
\ee 

Mapping grades $k$, $0 \leq k \leq 3$, to grades $3-k$ is known as \textit{duality} or \textit{Hodge duality}. Because of its usefulness and importance it gets the symbol $*$ as an upper index
\begin{gather}
  M^* := M i_3^{-1} = -m_0i_3 - m_1 e_{23} - m_2 e_{31} - m_3 e_{12} 
  \nonumber \\
  + m_{23} e_{1} + m_{31} e_{2}+ m_{12} e_{3} + m_{123}.
  \label{eq:Cl3duality}
\end{gather}
\textit{Duality} (\ref{eq:Cl3duality}) in $Cl(3,0)$ changes the outer product of vectors into the cross product
\begin{gather} 
  (e_1\wedge e_2)^* = e_1e_2 (-e_{123})= e_3 = e_1 \times e_2 , 
  \nonumber \\
  (e_2\wedge e_3)^* = e_1, \qquad 
  (e_3\wedge e_1)^* = e_2.
\end{gather}
Therefore we can use
\be 
  (a \wedge b)^* = a \times b, \qquad a \wedge b = i_3 (a \times b),
\ee 
to \textit{translate} well known results of standard 3D vector algebra into GA. Yet we emphasize, that the cross product $a \times b$ only exists in 3D, where a bivector has a unique (orthogonal) dual vector. The outer product has the advantage, that it can be \textit{universally} used in all dimensions, because geometrically speaking the parallelogram spanned by two linearly independent vectors is always well defined, independent of the dimension of the embedding space. 
 
Duality\footnote{
  Equations (\ref{eq:outcondual1})
apply in \textit{all} Clifford algebras, but (\ref{eq:outcondual2}) need some modification if the pseudoscalar is not central.} also relates outer product and contraction of any $A,B \in Cl(3,0)$
\begin{align}  
  A \wedge B^* &= (A\rfloor B)^*, \qquad
  A^* \wedge B = (A\lfloor B)^*,
  \label{eq:outcondual1}
  \\
  A \rfloor B^* &= 
  A^* \lfloor B = (A\wedge B)^*.
  \label{eq:outcondual2}
\end{align}
For vectors $a,b \in \R^3$ this can be rewritten as
\begin{align} 
  a \wedge b^* &= a^* \wedge b = (a \cdot b)^* , 
  \nonumber \\
  a \rfloor b^* &= a^* \lfloor b = (a \wedge b)^* .
\end{align}
Inserting the duality definition this gives
\begin{align} 
  a \wedge (i_3b) = (i_3a)\wedge b = i_3(a \cdot b),
  \nonumber \\
  a \rfloor (i_3b) = (i_3a) \lfloor b = i_3(a \wedge b),
\end{align}
which are handy relations in GA computations. 

Let us do two simple examples. First, for $a=b=e_1$. 
Then we have $i_3a=i_3b=e_{23}$, $a\cdot b = 1$, and get
\begin{gather}  
  a \wedge (i_3b) = e_1 \wedge e_{23}=i_3,
  \nonumber \\
  (i_3a)\wedge b = e_{23}\wedge e_1 = i_3,
  \quad
  i_3(a \cdot b) = i_3.
\end{gather} 
Second, for $a=e_1, b=e_2$. 
Then we have $i_3a=e_{23}$, $i_3b=e_{31}$,  
$a \wedge b = e_{12}$, and get
\begin{align}  
  a \rfloor (i_3b) &= e_1 \rfloor e_{31} = \langle e_1 e_{31} \rangle_1 = -e_3,
  \\
  (i_3a) \lfloor b &= e_{23} \lfloor e_2 = -e_3,
  \quad
  i_3(a \wedge b) = i_3 e_{12} = -e_3.
  \nonumber 
\end{align}

\subsection{The blade subspaces of $Cl(3,0)$}

A vector $a\in \R^3$ can be used to define a hyperplane $H(a)$ by 
\be
  H(a) = \{x\in \R^3 \mid x \cdot a= x \rfloor a= 0 \}, 
\ee
which is an example of an \textit{inner product null space} (IPNS) definition. 
It can also be used to define a line $L(a)$ by 
\be 
  L(a) = \{x\in \R^3 \mid x \wedge a = 0 \}, 
\ee 
which shows an \textit{outer product null space} (OPNS) definition. 
We can also use two linearly independent vectors to span a plane $P(a\wedge b)$. By that we mean, that
$x = \alpha a + \beta b, \alpha, \beta \in \R$, if and only if, $x\wedge a\wedge b = 0$. 
Let us briefly check that. 
Assume $x = \alpha a + \beta b, \alpha, \beta \in \R$, then
$x\wedge a\wedge b 
= (\alpha a + \beta b)\wedge a \wedge b
=  \alpha (a\wedge a) \wedge b - \beta a \wedge (b \wedge b) = 0$.
On the other hand, if $x$ has a component not in the plane spanned by $a,b$, then $x\wedge a \wedge b$ results in a non-zero trivector, the oriented parallelepiped volume spanned by $x,a,b$. 
Hence every plane spanned by two linearly independent vectors $a,b \in \R$ is in OPNS
\begin{gather}
  P = P(a\wedge b) = \{x \in \R^3 \mid x = \alpha a + \beta b, \,\forall \,\alpha, \beta \in \R\}
  \nonumber \\
  = \{x\in \R^3 \mid x \wedge a \wedge b = 0\}.
  \label{eq:planeab}
\end{gather}

Three linearly independent vectors $a,b,c \in \R^3$ always span $\R^3$, and their outer product with any fourth vector $x\in \R^3$ is zero by default (because in $Cl(3,0)$ no 4-vectors exist). Therefore formally $\R^3 = \{x\in \R^3 \mid x \wedge a \wedge b \wedge c = 0\}$.

These geometric facts motivate the notion of \textit{blade}. A blade $A$ is an element of $Cl(3,0)$, that can be written as the outer product $A=a_1\wedge \ldots \wedge a_k$ of $k$ linearly independent vectors $a_1, \ldots , a_k \in \R^3$, $0 \leq k \leq 3$. Another name for $k$-blade $A$ is \textit{simple} $k$-vector. As we have just seen, every $k$-blade fully defines and characterizes a $k$-dimensional vector subspace $V(A)\subset \R^3$. Such subspaces are therefore often simply called \textit{blade subspaces}. The GA of a blade subspace $Cl(V(A))$ is a natural subalgebra of $Cl(3,0)$. Examples are $Cl(V(e_1)) = Cl(1,0) \subset Cl(3,0)$, 
$Cl(V(e_{12})) = Cl(2,0) \subset Cl(3,0)$, and $Cl(V(e_{123})=\R^3) = Cl(3,0) \subseteq Cl(3,0)$, etc. 

The formulas for projection and rejection of vectors can now be generalized to the projection and rejection of subspaces, represented by their blades in their OPNS representations, where the grade of a blade shows the dimension of the blade subspace. The \textit{projection} of a blade $A$ onto a blade $B$ is given by
\be 
  P_B(A) = (A\rfloor B)B^{-1},
  \label{eq:proj}
\ee 
where the left contraction automatically takes care of the fact that it would not be meaningful to project a higher dimensional subspace $V(A)$ onto a lower dimensional subspace $V(B)$. 

This projection formula is a first example of how GA allows to \textit{solve} geometric problems \textit{by} means of elementary \textit{algebraic products}, rather than the often cumbersome conventional solution of systems of linear equations. 

The \textit{rejection} (orthogonal part) of a blade $A$ from a second blade $B$ is given by
\be 
  P^{\perp}_B(A) = (A\wedge B)B^{-1}.
  \label{eq:rej}  
\ee  

If we interpret a general multivector $A\in Cl(3,0)$ as a \textit{weighted sum of blades}, where the weights can also be interpreted as blade volumes, then by linearity the projection $P_B(A)$ [and rejection $P^{\perp}_B(A)$] formula applied to a multivector $A$ yields the weighted sum of blades, each projected onto [rejected from] the blade $B$. 

Let us consider the example of $A = e_1+e_3 + e_{12}+e_{23}$ and $B=e_{12}$, $B^{-1}=-e_{12}$. We get 
\begin{gather}
  P_B(A) = [(e_1+e_3 + e_{12}+e_{23})\rfloor e_{12}](-e_{12})
  \nonumber \\
  = [e_2 + 0 -1 +0](-e_{12}) = e_1 + e_{12},
\end{gather}
which is exactly what we expect, because the vector $e_1$ is in the plane $B=e_{12}$, $e_3$ is orthogonal to $e_{12}$ and is projected out, the bivector component $e_{12}$ is also in the plane $B=e_{12}$, but the bivector component $e_{23}$ is perpendicular to $e_{12}$ and is correctly projected out. Let us also compute the rejection
\begin{gather}
  P^{\perp}_B(A) = [(e_1+e_3 + e_{12}+e_{23})\wedge e_{12}](-e_{12})
  \nonumber \\
  = [0+e_3e_{12}+0+0](-e_{12}) = e_3,
\end{gather}
which is again correct, because only the vector $e_3$ is orthogonal to the plane $B=e_{12}$, whereas e.g. $e_{23}$ is not fully orthogonal\footnote{To include  $e_{23}$ one can simply compute $A-P_B(A) = e_3+e_{23}$.}, because it has the vector $e_2$ in common with $B=e_{12}$.

The general duality of the outer product and the left contraction in (\ref{eq:outcondual1}) and (\ref{eq:outcondual2}) has the remarkable consequence of \textit{duality of IPNS and OPNS} subspace representations, because for all blades $A\in Cl(3,0)$
\be 
  x \rfloor A = 0
  \,\,\, \Leftrightarrow \,\,\,
  x \wedge A^* = 0, \quad \forall \, x\in \R^3,
\ee 
i.e. because in general $(x \rfloor A)^*=x \wedge A^*$ for all $x,A\in Cl(3,0).$ And because of $(A^*)^* = A (-i_3)^2 = - A$, we also have 
\be 
  x \rfloor A^* = 0
  \,\,\, \Leftrightarrow \,\,\,
  x \wedge A = 0, \quad \forall \, x\in \R^3.
\ee 
So we can either represent a subspace in OPNS or IPNS by a blade $A \in Cl(3,0)$ as
\begin{align} 
  V_{\text{OPNS}}(A) 
  &= \{x \in \R^3 \mid x \wedge A = 0\}
  \\
  &= \{x \in \R^3 \mid x \rfloor A^* = 0\}
  = V_{\text{IPNS}}(A^*).
  \nonumber 
\end{align}
We also have the relationship
\begin{align}  
  V_{\text{IPNS}}(A) 
  &= \{x \in \R^3 \mid x \rfloor A = 0\}
  \\
  &= \{x \in \R^3 \mid x \wedge A^* = 0\}
  = V_{\text{OPNS}}(A^*).
  \nonumber 
\end{align}
In application contexts both representations are frequently used and it should always be taken care to \textit{clearly specify} if a blade is meant to represent a subspace in the OPNS or the IPNS representation. In GA software like CLUCalc or CLUViz \cite{CP:CLUCalc} this is done by an initial multivector interpretation command. 

The outer product itself \textit{joins} (or unifies) disjoint (orthogonal) blade subspaces, e.g., the two lines $L(a)$ and $L(b)$ are unified by the outer product to the plane $P(a\wedge b)$ in (\ref{eq:planeab}). The left contraction has the dual property to \textit{cut out} one subspace contained in a larger one and leave only the orthogonal complement. For example, if $A=e_1$, $B=e_{12}$, then $A\rfloor B = e_1\rfloor e_{12} = e_2$, the orthogonal complement of $V(A)$ in $V(B)$. In general in the OPNS, the \textit{orthogonal complement} of $V(A)$ in $V(B)$ is therefore give by $V(A\rfloor B)$.

Let us assume two blades $A,B \in Cl(3,0)$, which represent two blade subspaces $V(A), V(B) \subset \R^3$ with intersection $V(M) = V(A)\cap V(B)$, given by the common blade factor $M \in Cl(3,0)$ (called \textit{meet})
\begin{align}
  A &= A'M = A' \wedge M, 
  \nonumber \\
  B &= MB' = M\wedge B',
\end{align}
where $A',B'$ are the orthogonal complement blades of $M$ in $A, B$, respectively,
\begin{align}
  A' &= AM^{-1} = A \lfloor M^{-1}, 
  \nonumber \\ 
  B' &= M^{-1}B=M^{-1}\rfloor B,
\end{align}
where we freely use the fact that scalar factors like $M^{-2}=\pm|M|^{-2}\in \R$ of $M^{-1}=M\, M^{-2}$, are not relevant for determining a subspace (in both OPNS and IPNS), i.e. $V(M) = V(\lambda M), \, \forall\, \lambda \in \R\setminus \{0\}$.

The \textit{join} (union) $J$ of any two blade subspaces $V(A), V(B)$ can therefore be represented as
\be 
  J = A' \wedge M \wedge B' = A \wedge B' = A' \wedge B. 
\ee 
Inserting the above contraction formulas for $A', B'$ we get
\be 
  J = A \wedge (M^{-1}\rfloor B) = (A\lfloor M^{-1})\wedge B. 
\ee 

In turn we can use the join $J$ (or equivalently its inverse blade $J^{-1}$) to compute $M$. The argument is from set theory. If we cut out $B$ from $J$ (or $J^{-1}$), then only $A'$ (or $A'^{-1}$) will remain as the orthogonal complement of $B$ in $J$ (or $J^{-1}$). Cutting out $A'$ (or $A'^{-1}$) from $A$ itself leaves only the \textit{meet}\footnote{The symbol $\vee$ stems from Grassmann-Cayley algebra.} (intersection) $M$
\begin{align}
   M &= A \vee B = A'^{-1} \rfloor A = (B\rfloor J^{-1})\rfloor A 
   \nonumber \\
   &= B\lfloor (J^{-1}\lfloor A),
\end{align}
The right most form $M = B\lfloor (J^{-1}\lfloor A)$ is obtained by cutting $A$ out of from $J$ (or $J^{-1}$), and using the resulting $B'^{-1}=J^{-1}\lfloor A$, to get $M$ as orthogonal complement of $B'$ in $B$. 

For example, if we assume $A=e_{12}$, $B=e_{23}$, then the join is obviously $J=e_{123}=i_3$. This allows to compute the meet as
\begin{gather} 
  M = A \vee B 
  = (B\rfloor J^{-1})\rfloor A
  \nonumber \\
  = (e_{23}\rfloor(-e_{123}))\rfloor e_{12}
  = e_1\rfloor e_{12}=e_2,
\end{gather}  
which represents the line of intersection $L(e_2)$ of the two planes $P(e_{12})$ and $P(e_{23})$.

\subsection{Practically important relations in $Cl(3,0)$}

Let $A,B \in Cl(p,q,r)$. 
Some authors use angular brackets without index to indicate the scalar part of a multivector, e.g.
\be 
  \langle AB \rangle = \langle AB \rangle_0 = \langle BA \rangle_0. 
\ee 
For reversion (an involution) we have
\begin{align} 
  (AB)^{\sim} &= \widetilde{B} \widetilde{A}, \quad
  \tilde{a}=a, \quad
  \langle \widetilde{A} \rangle_0=\langle A \rangle_0=\langle A \rangle_0^{\sim}, 
  \nonumber \\
  \langle A \rangle_k^{\sim} &= (-1)^{k(k-1)/2}\langle A \rangle_k.
\end{align} 
Apart from reversion, there is the automorphism \textit{main involution} (or \textit{grade involution}) which maps all vectors $a \rightarrow \hat{a}=-a$, and therefore
\be 
  \widehat{a_1 \ldots a_s} = (-1)^s a_1 \ldots a_s,
\ee 
i.e. even blades are invariant under the grade involution, odd grade blades change sign. 
The composition of reversion and grade involution is called \textit{Clifford conjugation}
\begin{align} 
  \overline{A}&=\widetilde{(\hat{A})} = \widehat{(\tilde{A})}, \qquad
  \overline{a_1 \ldots a_s} = (-1)^s a_s \ldots a_1, 
  \nonumber \\
  \overline{\langle A \rangle_k} &= (-1)^{k(k+1)/2}\langle A \rangle_k.
\end{align} 
Applied to a multivector $M \in Cl(3,0)$ we would get
\begin{align} 
  \widetilde{M} 
  &= \langle M \rangle_0 + \langle M \rangle_1 -\langle M \rangle_2 -\langle M \rangle_3,
  \nonumber \\
  \widehat{M} 
  &= \langle M \rangle_0 - \langle M \rangle_1 +\langle M \rangle_2 -\langle M \rangle_3,
  \nonumber \\
  \overline{M} 
  &= \langle M \rangle_0 - \langle M \rangle_1 -\langle M \rangle_2 +\langle M \rangle_3.
\end{align} 
The following algebraic identities \cite{DH:NFCM} are frequently applied
\begin{align}
  &a\rfloor (b\wedge c)
  = (a\rfloor b)c -(a\rfloor c)b
  = (a\cdot b)c -(a\cdot c)b
  \nonumber \\
  &\hspace*{5mm}\stackrel{Cl(3,0)}{=} -a\times (b\times c),
  \\
  &a\rfloor (b\wedge c \wedge d)
  = (b\wedge c \wedge d)\lfloor a,
  \\
  &\hspace*{3mm}= (b\wedge c)(a\cdot d)-(b\wedge d)(c\cdot a)+(c\wedge d)(b\cdot a),
  \nonumber \\
  &a\wedge b \wedge c \wedge d \stackrel{Cl(3,0)}{=}0.
\end{align}

\section{Clifford's geometric algebra of spacetime (STA)}

The GA $Cl(1,3)$ of flat Minkowski \textit{spacetime}\footnote{Note that some authors prefer opposite signature $Cl(3,1).$} $\R^{1,3}$ \cite{DL:GAfPh}, has the $2^4=16$D basis ($e_0$ represents the time dimension)
\begin{gather}
   \{1, e_0, e_1, e_2, e_3, 
   \sigma_1=e_{10}, \sigma_2=e_{20}, \sigma_3=e_{30}, 
   \nonumber \\
   e_{23}, e_{31}, e_{12},
   e_{123}, e_{230}, e_{310}, e_{120},
   e_{0123}=I
   \}
\end{gather}
of one scalar, 4 vectors (basis of $\R^{1,3}$), 6 bivectors, 4 trivectors and one pseudoscalar. The 4 vectors fulfill 
  $e_{\mu} \rfloor e_{\nu} = \eta_{\mu,\nu}=\mathrm{diag}(+1,-1,-1,-1)_{\mu,\nu}\,,\, 0\leq \mu,\nu \leq 3$.

The commutator of two bivectors always gives a third bivector, e.g., 
\be 
  \frac{1}{2}(e_{10}e_{20}-e_{20}e_{10}) = -e_{12}, \,\ldots
\ee 
Exponentiating all bivectors gives a group of rotors, the \textit{Lorentz group}, its \textit{Lie algebra} is the commutator algebra of the six bivectors. The even grade subalgebra $Cl^+(1,3)$ with 8D basis $(l=1,2,3)$: $\{1, \{\sigma_l\}, \{I\sigma_l\}, I\}$ is isomorphic to $Cl(3,0)$. 

Assigning \textit{electromagnetic} field vectors $\vec{E}=E_1\sigma_1+E_2\sigma_2+E_3\sigma_3$,  $\vec{B}=B_1\sigma_1+B_2\sigma_2+B_3\sigma_3$, the Faraday bivector is $F=\vec{E}+I\vec{B}$. All four \textit{Maxwell equations} then \textit{unify} into one
\be 
  \nabla F = J,
\ee 
with vector derivative $\nabla = \sum_{\mu,\nu=0}^3\eta^{\mu,\nu} e_{\nu}\partial_{\mu}=\sum_{\mu=0}^3e^{\mu}\partial_{\mu}$. 

The \textit{Dirac} equation of \textit{quantum mechanics} for the electron becomes in STA simply
\be 
  \nabla \psi I \sigma_3 = m \psi e_0,
\ee  
where the spinor $\psi$ is an even grade multivector, and $m$ the mass of the electron. The electron spinor acts as a spacetime rotor resulting in observables, e.g., the current vector
\be 
  J = \psi e_0 \widetilde{\psi}.
\ee

\section{Conformal geometric algebra \label{sc:CGA}}

\subsection{Points, planes and motors in $Cl(4,1)$}

In order to \textit{linearize} translations as in (\ref{eq:lintrans}) we need the conformal
model~\cite{SL:Diss, LVA:MoebRn, HL:IAandGR} of Euclidean space (in $Cl(4,1)$), 
which adds to $\gvec{x}\in\R^3\subset\R^{4,1}$ two null-vector dimensions for 
the origin $\gvec{e}_0$ and infinity $\gvec{e}_{\infty}$, 
$\gvec{e}_0\rfloor \gvec{e}_{\infty}=-1$, to obtain by a non-linear embedding in $\R^{4,1}$ the homogeneous\footnote{Unique up to a nonzero scalar factor.} conformal point
\begin{align}
  &X = \gvec{x} + \frac{1}{2}\gvec{x}^2\gvec{e}_{\infty}+\gvec{e}_0, 
  \qquad
  \gvec{e}_0^2=\gvec{e}_{\infty}^2=X^2=0,
  \nonumber \\
  &X\rfloor \gvec{e}_{\infty}=-1.
\end{align}
The condition $X^2=0$ restricts the point manifold to the so-called null-cone of $\R^{4,1}$, similar to the light cone of special relativity. The second condition $X\rfloor \gvec{e}_{\infty}=\gvec{e}_0\rfloor\gvec{e}_{\infty}=-1$ further restricts all points to a hyperplane of $\R^{4,1}$, similar to points in projective geometry. The remaining point manifold is again 3D. See \cite{DL:GAfPh,DFM:GAfCS,HTBY:Carrier,HL:IAandGR,AM:LinGA,CP:CLUCalc,CP:GAappEng,GS:CMinGA} for details and illustrations. We can always move from the Euclidean representation $\gvec{x} \in \R^3$ to the conformal point $X \in \R^{4,1}$ by adding the $\gvec{e}_{0}$ and $\gvec{e}_{\infty}$ components, or by dropping them (projection (\ref{eq:proj}) with $i_3$, or rejection (\ref{eq:rej}) with $E=\gvec{e}_{\infty}\wedge\gvec{e}_{0}$). 

The contraction of two conformal points gives their \textit{Euclidean distance} and 
therefore a plane $m$ equidistant from two points 
$A=\gvec{a} + \frac{1}{2}\gvec{a}^2\gvec{e}_{\infty}+\gvec{e}_0,
B=\gvec{b} + \frac{1}{2}\gvec{b}^2\gvec{e}_{\infty}+\gvec{e}_0$ as
\begin{align}
  &X\rfloor A = -\frac{1}{2}(\gvec{x}-\gvec{a})^2 
  \\
  &\Rightarrow \,\,X\rfloor (A-B)=0,
  \,\,\,
  m=A-B \propto \gvec{n}+d\,\gvec{e}_{\infty},
  \nonumber 
\end{align}
where $\gvec{n}\in\R^3$ is a unit normal to the plane and $d$ its signed scalar distance 
from the origin, ``$\propto$" means proportional. Reflecting at two parallel planes $m,m^{\prime}$ with 
distance $\gvec{t}/2\in\R^3$ we get the \textit{transla}tion opera\textit{tor} 
(by $\gvec{t}\,$)
\begin{equation}
  X^{\prime} = m^{\prime}m\,X\,mm^{\prime} = T_{\gvec{t}}^{-1} X T_{\gvec{t}},
  \,\,\, T_{\gvec{t}}:= mm^{\prime}=1+\frac{1}{2}\gvec{t}\gvec{e}_{\infty}.
  \label{eq:lintrans}
\end{equation}
Reflection at two non-parallel planes $m,m^{\prime}$ yields the rotation around
the $m,m^{\prime}$-intersection by twice the angle subtended by $m,m^{\prime}$.

Group theoretically the conformal group $C(3)$ \cite{GS:CMinGA} is isomorphic to $O(4,1)$ \cite{Wiki:IOG} and the
Euclidean group $E(3)$ is the subgroup of $O(4,1)$  leaving infinity 
$\gvec{e}_{\infty}$ invariant. 
Now general translations \textit{and} rotations are both \textit{linearly} represented by geometric products of invertible vectors (called Clifford monomials, Lipschitz elements, versors, or simply motion operators = \textit{motors}). The commutator algebra of the bivectors of $Cl(4,1)$ constitutes the \textit{Lie algebra} of all \textit{conformal transformations} (exponentials of bivectors) of $\R^3$. Derivatives with respect to these bivector motion parameters allow \textit{motor optimization}, used in pose estimation, structure and motion estimation, and motion capture.

\subsection{Blade representations of geometric objects in $Cl(4,1)$}

Computer vision, computer graphics and robotics are interested in the intuitive \textit{geometric} OPNS \textit{meaning} of blades (outer products of conformal points $P_1, \ldots , P_4$) in conformal GA ($Pp$ = point pair)
\begin{gather}
  Pp = P_1 \wedge P_2, 
  \qquad
  Circle = P_1 \wedge P_2\wedge P_3, 
  \nonumber \\
  \textit{Sphere} = P_1 \wedge P_2\wedge P_3 \wedge P_4,
\end{gather}
homogeneously representing the \textit{point pair} $\{P_1, P_2\}$, the \textit{circle} through $\{P_1, P_2, P_3\}$, and the \textit{sphere} with surface points $\{P_1, P_2, P_3, P_4\}$. Abstractly these are 0D, 1D, 2D, and 3D spheres $S$ with center $\gvec{c}\in \R^3$, radius $r\geq 0$, and \textit{Euclidean} carrier blade directions $\mathbf{D}$: 1 for points $P$, distance $\gvec{d}\in \R^3$ of $P_2$ from $P_1$, circle plane bivector $\mathbf{i}_c$, and sphere volume trivector $i_s \propto e_{123}=i_3$, see \cite{HTBY:Carrier}. 
\begin{align} 
  S = \mathbf{D}\wedge \gvec{c} 
      &+[\frac{1}{2}(c^2+r^2)\mathbf{D} - \gvec{c} (\gvec{c}\rfloor \mathbf{D})]\gvec{e}_{\infty}
      \nonumber \\
      &+ \mathbf{D}\gvec{e}_0 + (\mathbf{D}\lfloor \gvec{c})E,
\end{align}
with the origin-infinity bivector $E=\gvec{e}_{\infty}\wedge \gvec{e}_0$.

A point $X \in \R^{4,1}$ is on the sphere $S$, if and only if $S\wedge X=0$. By duality this is equivalent to $X\rfloor S^*=0$, which is equivalent in Euclidean terms to $(\gvec{x}-\gvec{c})^2=r^2$, where $\gvec{x}$ is the Euclidean part of $X$, $\gvec{c} \in \R^3$ the center of $S$, and $r$ the radius of $S$.  

These objects are stretched to infinity (flattened) by wedging with $\gvec{e}_{\infty}$
\begin{align} 
  F &= S\wedge \gvec{e}_{\infty}
  = \mathbf{D}\wedge \gvec{c}\gvec{e}_{\infty}-\mathbf{D}E
  \nonumber \\
  &= \mathbf{D}\gvec{c}_{\perp}\gvec{e}_{\infty}-\mathbf{D}E,
\end{align}
where $\gvec{c}_{\perp}\in\R^3$ indicates the support vector (shortest distance from the origin): $\gvec{p}$ for finite-infinite point pair $P\wedge \gvec{e}_{\infty}$,
$\gvec{c}_{\perp}$ for the line $P_1 \wedge P_2\wedge \gvec{e}_{\infty}$ and the plane
$P_1 \wedge P_2\wedge P_3 \wedge \gvec{e}_{\infty}$, and $0$ for the 3D space $\R^3$ itself ($F=-i_sE\propto I_5$). All geometric entities can be \textit{extracted} easily (as derived and illustrated in \cite{HTBY:Carrier})
\begin{align} 
  \mathbf{D} &= -F\lfloor E, \qquad
  r^2 = \frac{S\widehat{S}}{\mathbf{D}^2}, 
  \nonumber \\
  \gvec{c} &= \mathbf{D}^{-1}[S\wedge(1+E)]\lfloor E . 
\end{align}
The conformal model is therefore also a complete super model for projective geometry. 

In the dual IPNS representation spheres (points for $r=0$) and planes become vectors in $\R^{4,1}$
\begin{align}
  \textit{Sphere}^* &= C - \frac{1}{2}r^2\gvec{e}_{\infty}, 
  \nonumber \\
  Plane^* &= \gvec{n}+d\gvec{e}_{\infty}
\end{align}
where $C$ is the conformal center point, $\gvec{n}:=\gvec{c}_{\perp}/|\gvec{c}_{\perp}|\propto \mathbf{i}_ci_3$, ``$\propto$" means proportional.

A point $X \in \R^{4,1}$ is on the plane $Plane^*$, if and only if  $X\rfloor Plane^*=0$, which is equivalent in Euclidean terms to $\gvec{x} \cdot \gvec{n} = d$, where $\gvec{x}$ is the Euclidean part of $X$, $\gvec{n} \in \R^3$ the oriented unit normal vector of $Plane^*$, and $d$ its signed scalar distance from the origin. 

Points are spheres with $r=0$, circles and lines are intersection bivectors 
\begin{align}
  Circle^* &= S^*_1\wedge S^*_2, 
  \nonumber \\
  Line^* &= Plane_1^*\wedge Plane_2^*,
\end{align}
of two sphere, and two plane vectors, respectively. Inversion at a sphere becomes
$X \rightarrow SXS = S^*XS^*$, inversion at two concentric spheres gives scaling. For example, the conformal center point $C$ of a sphere $S$ is
  $C = S \gvec{e}_{\infty} \,S$.

\section{Clifford analysis}

Multivector valued functions
  $f: \R^{p,q} \rightarrow Cl_{p,q}, \,\, p+q=n, $ have $2^n$ blade components
  $(f_A: \R^{p,q} \rightarrow \R)$
\begin{equation}\label{eq:MVfunc}
    f(\vect{x})  =  \sum_{A} f_{A}(\vect{x}) {\vect{e}}_{A}.
\end{equation}
We define the \textit{inner product} of two
$\R^n \rightarrow Cl_{n,0}$ functions  $f, g$ by
\begin{gather}
  (f,g) 
  = \int_{\R^n}f(\mbox{\boldmath $x$})
    \widetilde{g(\mbox{\boldmath $x$})}\;d^n\mbox{\boldmath $x$}
  \nonumber \\
  = \sum_{A,B}\mbox{\boldmath $e$}_A \widetilde{\mbox{\boldmath $e$}_B}
    \int_{\R^n}f_A (\mbox{\boldmath $x$})
    g_B (\mbox{\boldmath $x$})\;d^n\mbox{\boldmath $x$},
  \label{eq:mc2}
\end{gather}
and the \textit{norm} for functions in $L^2(\R^n;Cl_{n,0})
   = \{f: \R^n \rightarrow Cl_{n,0} \mid \|f\| < \infty \}$ as
\begin{equation}\label{eq:0mc2}
  \|f\|^2 
   = \left\langle ( f,f ) \right\rangle
   =
   \int_{\mathbb{R}^n} |f(\vect{x})|^2 d^n\mbox{\boldmath $x$}.
\end{equation}

The possibility to differentiate with respect to any multivector (representing a geometric object or a transformation operator) is essential for solving \textit{optimization} problems in GA. 

We can then define the \textit{vector differential} \cite{HS:CAtoGC,EH:VecDC} of $f$ for any constant $\vect{a} \in \mathbb{R}^{p,q}$ as
\begin{equation}\label{eqma1}
  \vect{a}\cdot\nabla f(\vect{x})=\lim_{\epsilon  \rightarrow 0}
  \frac{f(\vect{x}+\epsilon \vect{a})-f(\vect{x})}{\epsilon } , 
\end{equation}
where $\vect{a}\cdot\nabla$ is scalar.
The \textit{vector} derivative $\nabla$ can be expanded in terms of the basis vectors $\vect{e}_k$ as
\begin{align}  
  \nabla &= \nabla_{\scriptstyle \vect{x}} = \sum_{k=1}^{n} \vect{e}_k {\partial_k },
  \nonumber \\
  {\partial_k } &= \vect{e}_k \cdot\nabla = \frac{\partial}{\partial x_k}.
  \label{eq4}
\end{align}
Both $\vect{a}\cdot\nabla$ and $\nabla$ are \textit{coordinate independent}.
Replacing the vectors $\vect{a},\vect{x}$ by multivectors gives the \textit{multivector derivative} \cite{HS:CAtoGC,EH:MVDC}.

\textit{Examples.}
The five multivector functions 
\begin{align}
  f_1 &= \vect{x}, \quad
  f_2 = \vect{x}^2, \quad
  f_3 = |\vect{x}|,
  \nonumber \\
  f_4 &= \vect{x}\cdot \langle A \rangle_k, \qquad
    0\leq k \leq n,
  \\
  f_5 &= \log r, \quad 
  \vect{r}=\vect{x}-\vect{x}_0, \quad
  r = |\vect{r}|, 
  \nonumber 
\end{align}
have the following vector differentials \cite{EH:VecDC}
\begin{align}
  \vect{a}\cdot\nabla f_1 &= \vect{a}, \,\,
  \vect{a}\cdot\nabla f_2 = 2\vect{a}\cdot \vect{x}, \,\,
  \vect{a}\cdot\nabla f_3 = \frac{\vect{a}\cdot \vect{x}}{|\vect{x}|}, 
  \nonumber \\
  \vect{a}\cdot\nabla f_4 &= \vect{a}\cdot \langle A \rangle_k,\quad 
  \vect{a}\cdot\nabla f_5 = \frac{\vect{a}\cdot \vect{r}}{r^2}. 
\end{align}
This leads to the vector derivatives \cite{EH:VecDC}
\begin{align}
  \nabla f_1 &= \nabla\cdot \vect{x} = 3, \, (n=3), \,\,
  \nabla f_2 = 2 \vect{x}, \,\,
  \nabla f_3 = \vect{x}/|\vect{x}|, 
  \nonumber \\
  \nabla f_4 &= k \langle A \rangle_k, \qquad
  \nabla f_5 = \vect{r}^{-1}= \vect{r}/r^2.
\end{align}

We can compute the \textit{derivative from the differential} 
for $\nabla_{\scriptstyle \vect{a}}$: regard $\vect{x}=$ constant, $\vect{a}=$ variable, and compute
\begin{equation}
\nabla f = \nabla_{\scriptstyle \vect{a}} \,(\vect{a} \cdot \nabla f )
\end{equation}
The vector derivative obeys
\textit{sum rules} and \textit{product rules} (overdots indicate functions to be differentiated). But non-commutativity leads to modifications, because
$\dot{\nabla}f \dot{g} \neq f \dot{\nabla} \dot{g}$. The \textit{product rule} is
\begin{align}
  \nabla (fg)&=(\dot{\nabla}\dot{f} )g + \dot{\nabla}f \dot{g}
  \nonumber \\
  &=(\dot{\nabla}\dot{ f})g + \sum_{k=1}^n \vect{e}_k f ({\partial_k g}).
\end{align}
The \textit{chain rules} for the 
vector differential and the vector derivative 
of $f(\vect{x}) = g(\lambda (\vect{x})), \,\lambda (\vect{x})\in \R$, are
\begin{align}
  \vect{a} \cdot \nabla f 
  &= \{ \vect{a} \cdot \nabla \lambda ( \vect{x} ) \} 
    \frac{\partial g}{\partial{\lambda}},
  \nonumber \\
  \nabla f 
  &= (\nabla\lambda)\frac{\partial g}{\partial \lambda}.
\end{align}
\textit{Example:} 
For $\vect{a}=\vect{e}_k$ $(1\leq k \leq n)$ we have
\be 
  \vect{e}_k \cdot \nabla f
  = \partial_k f 
  = (\partial_k \lambda) \frac{\partial g}{\partial \lambda} .
\ee 
GA thus provides a new formalism for \textit{differentiation on vector manifolds} and for \textit{mappings} between surfaces, including \textit{conformal mappings} \cite{GS:CMinGA}.

Within Clifford analysis we can define \textit{quaternionic} and \textit{Clifford FTs} and \textit{wavelet transforms}, with applications in image and signal processing, and \textit{multivector wave packet analysis} \cite{EH:QFTgen}. We also obtain a single \textit{fundamental theorem of multivector calculus}, which \textit{unifies} a host of classical theorems of integration for path independence, Green's, Stokes', and Gauss' divergence theorems. \textit{Monogenic} functions $f$, $\nabla f = 0$, generalize complex analytic functions, and allow us to generalize Cauchy's famous integral theorem for analytic functions in the complex plane to $n$ dimensions \cite{SS:FThC}. 

\section*{Acknowledgments}

\textit{Now to the King eternal, immortal, invisible, the only God, be honor and glory for ever and ever. Amen.} [Bible, 1 Tim. 1:17] I do thank my family for their patient support, Y. Kuroe, T. Nitta and the anonymous reviewers. E.H. found several Wikipedia articles helpful.

\begin{biography}

\profile{m}{Eckhard Hitzer}{
  He received his M.S. from the Techn. Univ. of Munich, Germany, in 1992, and his Ph.D. from the Univ. of Konstanz, Germany, in 
  1996. In 1996 he joined the RIMS in Kyoto, Japan, for 2 years as a postdoc. 
  In 1998, he joined the Univ. of Fukui, where he is currently a part time lecturer in the Dep. of Appl. Phys.
  His research interests include Cliff. algebra and analysis, visualization, and neural computing.  
  He is a member of DPG, JPS, CAVi, CAIROS, ENNS, AMS, JSIAM, SICE and the Task Force on Complex-Valued NNs in the IEEE CIS NNs Techn. Comm.}  
\end{biography}

\end{document}